\documentclass[11pt]{amsart}

\usepackage{amsmath,amsthm}
\newtheorem{theorem}{Theorem}[section]
\newtheorem{lemma}[theorem]{Lemma}
\newtheorem{corollary}[theorem]{Corollary}
\newtheorem{prop}[theorem]{Proposition}

\theoremstyle{definition}

\theoremstyle{remark}
\newtheorem{remark}[theorem]{Remark}

\numberwithin{equation}{section}

\def\R{\mathbb{R}}
\def\e{\rm e}
\def\k{k{\rm !}}
\def\r{r{\rm !}}
\def\u{\rm u}

\def\N{\mathbb{N}}
\def\P{\mathbb{P}}
\def\rp{\mathbf{p}}
\def\y{\mathbf{y}}
\def\rq{\mathbf{q}}

\def\mP{\mathcal{P}}
\def\K{\mathbb{K}}
\def\Q{\mathbb{Q}}

\def\ring{\R[x_1,\ldots,x_n]}

\begin{document}

\title{A sum of squares approximation of nonnegative polynomials}

\author{Jean B. Lasserre}

\address{LAAS-CNRS\\
7 Avenue du Colonel Roche\\
31077 Toulouse C\'{e}dex 4, France}


\date{}


\begin{abstract}
We show that every real nonnegative polynomial $f$
can be approximated as closely as desired by a sequence of 
polynomials $\{f_\epsilon\}$ that are sums of squares. Each
$f_\epsilon$ has a simple et explicit form in terms of $f$ 
and $\epsilon$. A special representation
is also obtained for convex polynomials, 
nonnegative on a convex semi-algebraic set.
\end{abstract}

\maketitle

\section{Introduction}

The study of relationships between {\it nonnegative} and {\it sums of
squares} (s.o.s.) polynomials, initiated by Hilbert, is of real
practical importance in view of numerous potential applications,
notably in  polynomial programming. Indeed, checking whether a given polynomial is
nonnegative is a NP-hard problem whereas checking it is 
s.o.s. reduces to  solving a (convex) Semidefinite Programming (SDP) problem for
which efficient algorithms are now available. 
For instance, recent results in real algebraic geometry, most notably 
by Schm\"udgen \cite{schmudgen}, Putinar \cite{putinar}, 
Jacobi and Prestel \cite{jacobi}, have provided s.o.s. representations
of polynomials, positive on a compact semialgebraic set; the
interested reader is referred to Prestel and Delzell \cite{prestel}
and Scheiderer \cite{claude} for a nice account of such results. This in turn
has permitted to develop efficient SDP-relaxations in polynomial
optimization (see e.g. Lasserre \cite{lasserre1,lasserre2,lasserre3}, Parrilo
\cite{parrilo1,parrilo2}, Schweighofer \cite{markus}, and the many references therein).

So, back to a comparison between nonnegative and s.o.s. polynomials,
on the negative side, Blekherman \cite{blekherman} has shown that if
the degree is {\it fixed}, then the cone of nonnegative polynomials is much
{\it larger} than that of s.o.s. However, on the positive side, 
a denseness result \cite{berg}
states that the cone of s.o.s. polynomials is {\it dense}
in the space of polynomials that are nonnegative on $[-1,1]^n$ 
(for the norm $\Vert f\Vert_1=\sum_\alpha\vert f_\alpha\vert$ whenever
$f$ is written $\sum_\alpha f_\alpha x^\alpha$ in the usual canonical
basis); see e.g. Berg \cite[Theorem 5, p. 117]{berg}). 
\vspace{0.1cm}

{\bf Contribution.}
We show that {\it every nonnegative} polynomial $f$ is almost a s.o.s.,
namely we show that
$f$ can be approximated by a
sequence of s.o.s. polynomials
$\{f_\epsilon\}_\epsilon$, in the specific form
\begin{equation}
\label{specific}
f_\epsilon\,=\,f+\epsilon
\sum_{k=0}^{r_\epsilon}\sum_{j=1}^n \frac{x_j^{2k}}{\k},
\end{equation}
for some $r_\epsilon\in\N$, so that $\Vert f-f_\epsilon\Vert_1\to 0$ as $\epsilon\downarrow 0$.

This result is in the spirit of the previous denseness result.
However we here provide in (\ref{specific}) an {\it explicit} converging approximation with a very
specific (and simple) form; namely it suffices to slightly perturbate
$f$ by adding a small coefficient $\epsilon>0$ to each square monomial
$x^{2k}_i$ for all $i=1,\ldots,n$ and all $k=1,\ldots,r$, with $r$
sufficiently large.

To prove this result we combine 

- (generalized) {\bf Carleman}'s sufficient condition 
for a moment sequence $\y=\{y_\alpha\}$ to have a {\it
representing measure} $\mu$ (i.e., such that $y_\alpha=\int x^\alpha
d\mu$ for all $\alpha\in\N^n$), and

- a {\bf duality} result from convex optimization.

As a consequence, we may thus define a procedure to
approximate the global minimum of a polynomial $f$.
It consists in solving a sequence of SDP-relaxations which are simpler and easier to solve
than those defined in Lasserre \cite{lasserre1}.
Finally, we also consider the case where $f$ is a {\it convex} polynomial,
nonnegative on a convex semi-algebraic set $\K$ defined by (concave polynomial)
inequalities $g_j\geq 0$. We show that the approximation $f_\epsilon$
of $f$, defined in (\ref{specific}), has
a certificate of positivity on $\K$ (or a representation) similar to
Putinar's s.o.s. representation \cite{putinar}, but in which the
s.o.s. polynomial coefficients of the $g_j$'s now become simple
nonnegative {\it scalars}, the Lagrange multipliers of a related convex
optimization problem.

\section{Notation and definitions}
\label{notation}

For a real symmetric matrix $A$, the notation $A\succeq0$
(resp. $A\succ0$) stands for
$A$ positive semidefinite (resp. positive definite).
The sup-norm $\sup_j\vert x_j\vert$ of a 
vector $x\in\R^n$, is denoted by $\Vert x\Vert_\infty$.
Let $\ring$ be the ring of real polynomials, and let

\begin{equation}
\label{aa4}
v_r(x):=(1,x_1,x_2,\dots x_n,x_1^2,x_1x_2,\dots,x_1x_n,x_2^2,x_2x_3,\dots,x_n^r)
\end{equation}
be the canonical basis for the $\R$-vector space 
$\mathcal{A}_r$ of real polynomials of degree at most $r$,
and let $s(r)$ be its dimension.
Similarly, $v_\infty(x)$ denotes the canonical basis of $\ring$ as a $\R$-vector
space, denoted $\mathcal{A}$. So a vector in $\mathcal{A}$ has always
{\it finitely} many zeros.

Therefore, a polynomial $p\in\mathcal{A}_r$ is written
\[x\mapsto p(x)\,=\,\sum_{\alpha}p_\alpha
x^\alpha\,=\,\langle \rp,v_r(x)\rangle,\hspace{1cm}x\in\R^n,\]
(where $x^{\alpha}=x_1^{\alpha_1} x_2^{\alpha_2} \dots
x_n^{\alpha_n}$) for some vector 
$\rp=\{ p_\alpha\}\in\R^{s(r)}$, the vector of coefficients of
$p$ in the basis (\ref{aa4}). 

Extending $\rp$ with zeros, we can also
consider $\rp$ as a vector indexed in the basis $v_\infty(x)$
(i.e. $\rp\in\mathcal{A}$). 
If we equip $\mathcal{A}$ with the usual scalar product $\langle
.,.\rangle$ of vectors, then 
for every $p\in\mathcal{A}$,
\[p(x)\,=\,\sum_{\alpha>\in\N^n}p_\alpha
x^\alpha\,=\,\langle \rp,v_\infty(x)\rangle,\hspace{1cm}x\in\R^n.\]

Given a sequence $\y=\{y_\alpha\}$ indexed in the basis $v_\infty(x)$,
let $L_\y:\mathcal{A}\to\,\R$ be the linear functional
\[p\mapsto L_\y(p)\,:=\,\sum_{\alpha\in\N^n}\,p_\alpha
y_\alpha\,=\,\langle \rp ,\y\rangle.\]
Given a sequence $\y=\{y_\alpha\}$ indexed in the basis $v_\infty(x)$, the {\it moment}
matrix $M_r(\y)\in\R^{s(r)\times s(r)}$ with
rows and columns indexed in the basis $v_r(x)$ in (\ref{aa4}), satisfies
\[\left[M_r(\y)(1,j)\,=\,y_\alpha \:\mbox{ and
}\:M_r(y)(i,1)\,=\,y_\beta\right]\,\Rightarrow\,
M_r(y)(i,j)\,=\,y_{\alpha+\beta}.\]
For instance, with $n=2$, 
\[M_2(\y)\,=\,\left[\begin{array}{cccccc}
y_{00}&y_{10}&y_{01} &y_{20}&y_{11}&y_{02}\\
y_{10}&y_{20}&y_{11}&y_{30}&y_{21}&y_{12}\\
y_{01}&y_{11}&y_{02}&y_{21}&y_{12}&y_{03}\\
y_{20}&y_{30}&y_{21}&y_{40}&y_{31}&y_{22}\\
y_{11}&y_{21}&y_{12}&y_{31}&y_{22}&y_{13}\\
y_{02}&y_{12}&y_{03}&y_{22}&y_{13}&y_{04}
\end{array}\right].\]
A sequence $\y=\{y_\alpha\}$ has a {\it representing} measure $\mu_\y$
if
\begin{equation}
\label{momentseq}
y_\alpha\,=\,\int_{\R^n} x^\alpha\,d\mu_\y,\qquad \forall\,\alpha\in\N^n.
\end{equation} 
In this case one also says that $\y$ is a {\it moment sequence}. In addition, if
$\mu_\y$ is unique then $\y$ is said to be a {\it determinate} moment sequence.

The matrix $M_r(\y)$ defines a bilinear form $\langle .,.\rangle_\y$ 
on $\mathcal{A}_r$, by
\[\langle q,p\rangle_\y\,:=\,\langle
\rq,M_r(\y)\rp\rangle\,=\,L_\y(qp),\hspace{0.3cm}q,p\in\mathcal{A}_r,\]
and if $\y$ has a {\it representing} measure $\mu_\y$ then
\begin{equation}
\label{moment1}
\langle \rq,M_r(\y)\rq\rangle\,=\,\int_{\R ^n} q(x)^2\,\mu_\y(dx)\,\geq\,0,
\end{equation}
so that $M_r(\y)\succeq 0$. 

Next, given a sequence $\y=\{y_\alpha \}$ indexed in the basis $v_\infty(x)$, 
let $y^{(i)}_{2k}:=L_\y(x_i ^{2k})$ 
for every $i=1,\ldots ,n$ and every $k\in\N$. That is, $y^{(i)}_{2k}$
denotes the element in the sequence $\y$, corresponding to the monomial
$x_i ^{2k}$.

Of course not every sequence $\y=\{y_\alpha\}$ has a representing
measure $\mu_\y$ as in (\ref{momentseq}). However, there exists a
{\it sufficient} condition to ensure that it is the case. The
following result stated in Berg \cite[Theorem 5, p. 117]{berg} is from
Nussbaum \cite{nussbaum}, and is re-stated here, with our notation.

\begin{theorem}
\label{suff}
Let $\y=\{y_\alpha\}$ be an infinite sequence such that
$M_r(\y)\succeq 0$ for all $r=0,1,\ldots$. If
\begin{equation}
 \label{suff-1}
\sum_{k=0}^\infty (\y^{(i)}_{2k})^{-1/2k}\,=\,\infty,\qquad i=1,\ldots,n,
 \end{equation}
then  $\y$ is a determinate moment sequence.
\end{theorem}

The condition (\ref{suff-1}) in Theorem \ref{suff} is called
{\it Carleman's condition} as it extends to the multivariate case 
the original Carleman's sufficient condition given for the univariate case.

\section{Preliminaries}

Let $B_M$ be the closed ball
\begin{equation}
\label{ball}
B_M\,=\,\{x\in\R^n\,\vert\quad \Vert x\Vert_\infty \,\leq \,M\}.
\end{equation}

\begin{prop}
\label{prop0}
Let $f\in\ring$ be such that $-\infty<f^*:=\inf_{x}f(x)$. 
Then, for every $\epsilon >0$ there is
some $M_\epsilon\in\N$ such that 
\[f^*_M\,:=\,\inf_{x\in B_M}\,f(x)
\,<\,f^*+\epsilon,\qquad \forall M\,\geq \,M_\epsilon.\]
Equivalently, $f^*_M\downarrow f^*$ as $M\to\infty$.
\end{prop}
\begin{proof}
Suppose it is false. That is, there is some $\epsilon_0>0$ and an
infinite sequence sequence $\{M_k\}\subset\N$, with $M_k\to\infty$, such that 
$f^*_{M_k}\geq f^*+\epsilon_0$ for all $k$. But let $x_0\in\R^n$ be such
that $f(x_0)<f^*+\epsilon_0$. With any $M_k\geq \Vert x_0\Vert_\infty$,
one obtains the contradiction $f^*+\epsilon_0\leq f^*_{M_k}\leq f(x_0)<f^*-\epsilon_0$.
\end{proof}

To prove our main result (Theorem \ref{uncons} below),
we first introduce the following related optimization problems.
\begin{equation}
\label{a0}
\P:\qquad f^*\,:=\,\inf_{x\in\R^n}\,f(x),
\end{equation}
and for $0<M\in\N$, 
\begin{equation}
\label{a1}
\mathcal{P}_M: \:\displaystyle{\inf_{\mu\in\mathcal{P}(\R^n)}}
\left\{\int f\,d\mu \,\vert\quad 
\int\sum_{i=1}^n \e ^{x^2_i}\,d\mu\leq n\e ^{M^2}\right\},
\end{equation}
where $\mathcal{P}(\R^n)$ is the space of probability measures on $\R^n$.
The respective optimal values of $\P$ and $\mathcal{P}_M$ are denoted $\inf\P=f^*$ and
$\inf\mathcal{P}_M$, or $\min\P$ and 
$\min\mathcal{P}_M$ if the minimum is attained.

\begin{prop}
\label{prop1}
Let $f\in\ring$ be such that $-\infty< f^*:=\inf_x f(x)$, and
consider the two optimization problems $\P$ and $\mP_M$ defined in
(\ref{a0}) and (\ref{a1}) respectively.
Then, $\inf\mathcal{P}_M\downarrow f^*$ as $M\to\infty$. If $f$ has a
global minimizer $x^*\in\R^n$, then $\min\mathcal{P}_M=f^*$ whenever
$M\geq \Vert x^*\Vert_\infty$.
\end{prop}

\begin{proof}
Let $\mu\in\mP(\R^n)$ be admissible for $\mP_M$.
As $f\geq f^*$ on $\R^n$ then it follows immediately that
$\int fd\mu\geq f^*$, and so, $\inf\mathcal{P}_M\geq f^*$
for all $M$. 

As $B_M$ is closed and bounded, it is compact and so, with $f^*_M$ as
in Proposition \ref{prop0}, there is some
$\hat{x}\in B_M$ such that $f(\hat{x})=f^*_M$. In addition
let $\mu\in\mP(\R^n)$ be the Dirac probability measure at the
point $\hat{x}$. As $\Vert\hat{x}\Vert_\infty\leq M$,
\[\int\,\sum_{i=1}^n \e ^{x^2_i}\,d\mu\,=\,
\sum_{i=1}^n \e ^{(\hat{x}_i)^2}\,\leq\, n\e ^{M^2},\]
so that $\mu$ is an admissible solution of $\mathcal{P}_M$
with value $\int f\,d\mu=f(\hat{x})=f^*_M$, which proves that $\inf\mathcal{P}_M
\leq f^*_M$. This latter fact, combined with 
Proposition \ref{prop0} and with $f^*\leq\inf\mathcal{P}_M$, implies
$\inf\mathcal{P}_M\downarrow f^*$ as $M\to\infty$, the desired
result. The final statement is immediate by taking as feasible
solution for $\mathcal{P}_M$, the Dirac probability measure at the point $x^*\in B_M$
(with $M\geq \Vert x^*\Vert_\infty$). As its value is now $f^*$, it is
also optimal, and so, $\mathcal{P}_M$ is solvable with optimal value $\min\mathcal{P}_M=f^*$.
\end{proof}

Proposition \ref{prop1} provides a rationale for introducing the following Semidefinite
Programming (SDP) problems.
Let $2r_f$ be the degree of $f$ and for every $r_f\leq r\in\N$,
consider the SDP problem
\begin{equation}
\label{primalagain}
\Q_r\left\{\begin{array}{llcl}
&\displaystyle{\min_\y\:L_\y(f)\,(=\sum_{\alpha}f_\alpha y_\alpha})&&\\
{\rm s.t.}&M_{r}(\y)&\succeq &0\\
&\displaystyle{\sum_{k=0}^ r\sum_{i=1}^n \frac{y^{(i)}_{2k}}{\k}}&\leq &n\e ^{M^2},\\
&y_{0}&=&1,
\end{array}\right.
\end{equation}
and its associated {\it dual} SDP problem

\begin{equation}
\label{dualagain}
\Q^*_r\left\{\begin{array}{llc}
&\displaystyle{\max_{\lambda\geq 0,\gamma,q} \gamma -n\e ^{M^2}\lambda} &\\
{\rm s.t.}&f-\gamma&= \displaystyle{q-\lambda \sum_{k=0}^r\sum_{j=1}^n \frac{x_j^{2k}}{\k}}\\
&q &\mbox{s.o.s. of degree }\leq 2r,
\end{array}\right.
\end{equation}
with respective optimal values $\inf\Q_r$ and $\sup\Q^*_r$
(or $\min\Q_r$ and $\max\Q^*_r$ if the optimum is attained, in which case the problems
are said to be solvable). For more details on SDP theory, the interested
reader is referred to the survey paper \cite{boyd}.

The SDP problem $\Q_r$ is a relaxation of $\mathcal{P}_M$, and we next
show that in fact 

- $\Q_r$ is solvable for all $r\geq r_0$, 

- its optimal value $\min\Q_r\to \inf\mathcal{P}_M$ as $r\to\infty$, and

- $\Q^*_r$ is also solvable with same optimal value as $\Q_r$, for every
$r\geq r_f$.

This latter fact will be crucial to prove our main result in
the next section.
Let $l_\infty$ (resp. $l_1$) be the Banach space of bounded (resp. summable)
infinite sequences with the sup-norm  (resp. the $l_1$-norm).

\begin{theorem}
\label{th1}
Let $f\in\ring$ be of degree $2r_f$, with global
minimum $f^*>-\infty$, and let $M>0$ be fixed. Then :

{\rm (i)} For every $r\geq r_f$, $\Q_r$ is solvable, and
$\min\Q_r\uparrow \inf\mathcal{P}_M$ as $r\to\infty$.

{\rm (ii)} Let $\y^{(r)}=\{y^{(r)}_\alpha\}$ be an optimal solution
of $\Q_r$ and complete $\y^{(r)}$ with zeros to make it an element of
$l_\infty$. Every (pointwise) accumulation point $\y ^*$
of the sequence $\{\y^{(r)}\}_{r\in\N}$ is a determinate moment
sequence, that is,
\begin{equation}
y^*_\alpha\,=\,\int_{\R^n}x^\alpha\,d\mu ^*,\qquad \alpha\in\N^n,
\end{equation}
for a unique probability measure $\mu ^*$, and $\mu ^*$ is an optimal solution of
$\mathcal{P}_M$.

{\rm (iii)} For every $r\geq r_f$, $\max\Q^*_r=\min\Q_r$.
\end{theorem}
For a proof see \S\ref{proofth1}.

So, one can approximate the optimal value $f ^*$ of $\P$ as closely as
desired, by solving SDP-relaxations $\{\Q_r\}$ for sufficiently large
values of $r$ and $M$. Indeed, $f^*\leq \inf\mathcal{P}_M\leq f^*_M$,
with $f^*_M$ as in Proposition \ref{prop0}. 
Therefore, let $\epsilon >0$ be fixed, arbitrary. By
Proposition \ref{prop1}, we have 
$f^*\leq \inf\mathcal{P}_M\leq f^*+\epsilon$ provided that
$M$ is sufficiently large. Next, by Theorem \ref{th1}(i), one has 
$\inf\Q_r\geq\inf\mathcal{P}_M-\epsilon$ provided that $r$ is
sufficiently large, in which case, we finally have 
$f^*-\epsilon\leq \inf\Q_r\leq f ^*+\epsilon$.

Notice that the SDP-relaxation $\Q_r$ in (\ref{primalagain}) is simpler than the one defined
in Lasserre \cite{lasserre1}. Both have the same variables
$\y\in\R^{s(r)}$, but the former has {\it one} SDP constraint
$M_r(\y)\succeq 0$ and one scalar inequality (as one substitutes $y_0$
with $1$) whereas the latter has the same SDP constraint
$M_r(\y)\succeq 0$ and one additional SDP constraint $M_{r-1}(\theta
\y)\succeq0$ for the localizing matrix associated with the polynomial
$x\mapsto \theta(x)=M^2-\Vert x\Vert^2$. This results in a significant simplification.

\section{Sum of squares approximation}

Let $\mathcal{A}$ be equipped with the norm
\[f\,\mapsto \,\Vert f\Vert_1 \,:=\,\sum_{\alpha\in\N^n}\vert
f_\alpha\vert,\qquad f\in\mathcal{A}.\]

\begin{theorem}
\label{uncons}
Let $f\in\ring$ be nonnegative with global minimum $f^*$, 
that is,
\[0\,\leq\,f^*\,\leq\, f(x),\qquad x\in\R^n.\]

{\rm (i)} There is some $r_0\in\N,\lambda_0\geq0$
such that, for all $r\geq r_0$ and $\lambda\geq\lambda_0$,
\begin{equation}
\label{uncons-1}
f+\lambda\sum_{k=0}^r\sum_{j=1}^n \frac{x_j^{2k}}{\k}
\qquad\mbox{is a sum of squares.}
\end{equation}

{\rm (ii)} For every $\epsilon >0$, there is some $r_\epsilon\in\N$
such that, 
\begin{equation}
\label{uncons-2}
f_\epsilon:=f+\epsilon\sum_{k=0}^{r_\epsilon}\sum_{j=1}^n \frac{x_j^{2k}}{\k}
\qquad\mbox{is a sum of squares.}
\end{equation}
Hence, $\Vert f-f_\epsilon\Vert_1 \to 0$ as $\epsilon\downarrow 0$.
\end{theorem}
For a proof see \S\ref{proofuncons}.

\begin{remark}
\label{rem1}
Theorem \ref{uncons}(ii) is a {\it denseness} result in the spirit
of Theorem 5 in Berg \cite[p. 117]{berg}
which states that the cone of 
s.o.s. polynomials is dense (also for the norm $\Vert f\Vert_1$)
in the cone of polynomials that are nonnegative on $[-1,1]^n$.
However, notice that Theorem \ref{uncons}(ii) provides an {\it explicit} converging
sequence $\{f_\epsilon\}$ with a simple and very specific form.
\end{remark}

We next consider the case of a convex polynomial, nonnegative on a
convex semi-algebraic set. Given $\{g_j\}_{j=1}^m \subset\ring$, let $\K\subset\R^n$ be the
semi-algebraic set
\begin{equation}
\label{setk}
\K:=\{x\in\R^n\,\vert\quad g_j(x)\geq0,\quad j=1,\ldots,m\}.
\end{equation}

\begin{corollary}
\label{coro1}
Let $\K$ be as in (\ref{setk}), where all the $g_j$'s are concave, and
assume that Slater's condition holds, i.e., there exists $x_0\in\K$
such that $g_j(x_0)>0$ for all $j=1,\ldots,m$.

Let $f\in\ring$ be convex, nonnegative on $\K$, and with a minimizer
on $\K$, that is, $f(x^*)\leq f(x)$ for all $x\in\K$, for some $x ^*\in\K$. Then there exists
a nonnegative vector $\lambda\in\R^m$ such that for every
$\epsilon>0$, there is some $r_\epsilon\in\N$ for which
\begin{equation}
\label{coro-1}
f+\epsilon\sum_{k=0}^{r_\epsilon}\sum_{i=1}^n\frac{x_i^{2k}}{k{\rm
!}}\,=\,f_0+\sum_{j=1}^m\lambda_j\,g_j,
\end{equation}
with $f_0\in\ring$ being a sum of squares. (Therefore, the degree of
$f_0$ is less than $\max[2r_\epsilon, {\rm deg}\,f,{\rm deg}\,g_1,\ldots,{\rm deg}\,g_m]$.)
\end{corollary}

\begin{proof}
Consider the convex optimization problem
$f^*:=\min\{f(x)\:\vert\: x\in\K\}$.
As $f$ is convex, $\K$ is a convex set and Slater's condition holds,
the Karush-Kuhn-Tucker optimality
condition holds. That is, there exists a {\it nonnegative} vector
$\lambda\in\R^m$ of Lagrange-KKT multipliers, such that
\[\nabla f(x^*)\,=\,\sum_{j=1}^m\lambda_j\nabla g_j(x^*);\quad
\lambda_jg_j(x ^*)\,=\,0,\:j=1,\ldots,m.\]
(See e.g. Rockafellar \cite{convex}.)
In other words, $x ^*$ is also a (global) minimizer of the convex Lagrangian
$L:=f-\sum_{j=1}^m\lambda_jg_j$. Then $f^*=f(x^*)=L(x^*)$ 
is the (global) minimum of $f$ on $\K$,
as well as the global minimum of $L$ on $\R^n$, i.e.,
\begin{equation}
\label{lagrange}
f-\sum_{j=1}^m\lambda_jg_j -f^*\geq 0,\qquad x\in\R^n.
\end{equation}
As $f\geq0$ on $\K$, $f^*\geq0$, and so
$L\geq0$ on $\R^n$. Then
(\ref{coro-1}) follows from Theorem \ref{uncons}(ii), applied to the
polynomial $L$.
\end{proof}

When $\K$ is {\it compact} (and so, $f$ has necessarily a minimizer $x^*\in\K$),
one may compare Corollary \ref{coro1} with
Putinar's representation \cite{putinar}
of polynomials, positive on $\K$.
When $f$ is nonnegative on $\K$ (compact), and with
\begin{equation}
\label{approx}
f_\epsilon\,:=\,f+\epsilon\sum_{k=0}^{r_\epsilon}\sum_{i=1}^n\frac{x_i^{2k}}{k{\rm !}},
\end{equation}
one may rewrite (\ref{coro-1}) as
\begin{equation}
\label{rew}
f_\epsilon\,=\,f_0+\sum_{j=1}^m\lambda_j\,g_j,
%-\epsilon
\end{equation}
which is indeed a {\it certificate} of positivitity of $f_\epsilon$ on
$\K$. In fact, as $f_\epsilon>0$ on $\K$, (\ref{rew}) can be seen as a special form
of Putinar's s.o.s. representation, namely
\begin{equation}
\label{final}
f_\epsilon\,=\,q_0+\sum_{j=1}^mq_j\,g_j,\qquad \mbox{with $q_0,\ldots,q_m$
s.o.s.}
\end{equation}
(which holds under an additional assumption on the $g_j$'s).
So, in the convex compact case, and under Slater's condition,
Corollary \ref{coro1} states that if $f\geq0$ on $\K$, then its 
approximation $f_\epsilon$ in (\ref{approx}),
%\[f_\epsilon\,=\,f+\epsilon\sum_{k=0}^{r_\epsilon}\sum_{i=1}^n\frac{x^{2k}_i}{k{\rm!}},\]
has the simplified Putinar representation (\ref{rew}), in which the s.o.s. coefficients
$\{q_j\}$ of the $g_j$'s in (\ref{final}), become now simple nonnegative
{\it scalars} in (\ref{rew}), namely, the Lagrange-KKT multipliers $\{\lambda_j\}$.

\section{Proofs}
\label{proofs}

\subsection{Proof of Theorem \ref{th1}}
\label{proofth1}
We will prove (i) and (ii) together.
We first prove that $\Q_r$ is solvable. This is because the feasible
set (which is closed) is compact. Indeed, the constraint
\[\sum_{k=0}^ r\sum_{i=1}^n y^{(i)}_{2k}/k{\rm !}\,\leq \,n\e
^{M^2}\]
implies that every diagonal element $y^{(i)}_{2k}$ of
of $M_r(\y)$ is bounded by $\tau_r:=n\r\e ^{M^2}$. By Lemma \ref{lemma1}, this 
in turn implies that its diagonal elements (i.e.,
$y_{2\alpha}$, with $\vert\alpha\vert\leq r$) are {\it all} bounded by $\tau_r$.

This latter fact and again $M_r(\y)\succeq 0$, also imply that in fact {\it every}
element of $M_r(\y)$ is bounded by $\tau_r$,
that is, $\vert y_\alpha\vert\leq \tau_r$
for all $\vert\alpha\vert\leq 2r$. Indeed, for a
a symmetric matrix $A\succeq 0$, every non diagonal element $A_{ij}$ satisfies
$A_{ij}^2\leq A_{ii}A_{jj}$ so that $\vert A_{ij}\vert\leq \max_iA_{ii}$.

Therefore the set of feasible solutions of $\Q_r$ is a
closed bounded subset of $\R^{s(r)}$, hence compact. As $L_\y(f)$ is
linear in $\y$, the infimum is attained at some feasible point.
Thus, for all $r\geq r_f$, $\Q_r$ is solvable with optimal value
$\min\Q_r\leq \inf\mathcal{P}_M$. The latter inequality is because
the moment sequence $\y$ associated with an
an arbitrary feasible solution $\mu$ of $\mathcal{P}_M$, is obviously feasible
for $\Q_r$, and with value $L_{\y}(f)=\int fd\mu$.

Next, as the sequence
$\{\min\Q_r\}_r$ is obviously monotone non decreasing, one has
$\min\Q_r\uparrow \rho ^*\leq \inf\mathcal{P}_M$, as $r\to\infty$. 
We have seen that every entry of $M_r(\y)$ is bounded by $\tau_r$,
and this bound holds for all $r\geq r_f$.
Moreover, $M_r(\y)$ is also a (north-west corner)
submatrix of $M_s(\y)$ for every $s>r$. Indeed, whenever $s>r$, one may write
\[M_s(\y)\,=\,\left[\begin{array}{ccc}
M_r(\y)&\vert& B\\
-&\vert &-\\
B' &\vert & C\end{array}\right]\]
for some appropriate matrices $B$ and $C$. Therefore,
for the same reasons, any
feasible solution $\y$ of $\Q_s(\y)$ satisfies 
$\vert y_\alpha\vert\leq \tau_r$, 
for all $\alpha\in\N^n$ such that $\vert\alpha\vert\leq 2r$. Therefore, for every $s\in\N$,  and
every feasible solution $\y$ of $\Q_s$, we have
\[\vert y_\alpha\vert\leq \tau_r,\quad\forall \alpha\in\N^n,\,
2r-1\leq\vert \alpha\vert\leq 2r,\quad r=1,\ldots,s.\]
Thus, given $\y=\{y_\alpha\}$, denote by
$\hat{\y}=\{\hat{\y}_\alpha\}$ the new sequence obtained from $\y$ by
the scaling
\[\hat{\y}_\alpha\,:=\,\y_\alpha/\tau_r\qquad \forall
\alpha\in\N^n,\:2r-1\,\leq \vert\alpha\vert\leq 2r,\quad r=1,2,\ldots\]
So let $\y^{(r)}=\{y^{(r)}_\alpha\}$ be an optimal solution
of $\Q_r$ and complete $\y^{(r)}$ with zeros to make it an element of
$l_\infty$. 
Hence, all the elements $\hat{\y}^{(r)}$ are in the unit ball $B_1$ of
$l_\infty$, defined by 
\[B_1\,=\,\{\y=\{y_\alpha\}\in l_\infty\,\vert\quad \Vert
\y\Vert_\infty\leq 1\}.\]
By the Banach-Alaoglu Theorem, this ball is sequentially compact in
the $\sigma(l_\infty,l_1)$ (weak*) topology of $l_\infty$ (see e.g. Ash
\cite{ash}).
In other words, there exists an element $\hat{\y}^*\in B_1$ and a
subsequence $\{r_k\}\subset\N$, such that $\hat{\y}^{(r_k)}\to \hat{\y}^*$ for the
weak* topology of $l_\infty$, that is, for all $\u\in l_1$,
\begin{equation}
\label{weak*}
\langle \hat{\y} ^{(r_k)},\u\rangle\,\to\,
\langle \hat{\y} ^*,\u\rangle,\qquad \mbox{as }k\to\infty.
\end{equation}
In particular, {\it pointwise} convergence holds, that is, for all $\alpha\in\N^n$,
\[\hat{y}^{(r_k)}_\alpha\,\to\,\hat{y}^*_\alpha,\qquad \mbox{as
}k\to\infty,\]
and so, defining $\y^*$ from $\hat{\y}^*$ by
\[\y^*_\alpha\,=\,\tau_r\,\hat{\y}^*_\alpha,
\qquad \forall
\alpha\in\N^n,\:2r-1\,\leq \vert\alpha\vert\leq 2r,\quad
r=1,2,\ldots\]
one also obtains the pointwise convergence
\begin{equation}
\label{pointwise}
\mbox{for all}\quad\alpha\in\N^n, \quad y^{(r_k)}_\alpha\,\to\,y^*_\alpha,\qquad \mbox{as
}k\to\infty.
\end{equation}
We next prove that $\y ^*$ is the moment sequence of an optimal
solution $\mu ^*$ of problem
$\mathcal{P}_M$. 
From the pointwise convergence (\ref{pointwise}), we immediately
get $M_r(\y^*)\succeq 0$ for all $r\geq r_f$, 
because $M_r(\y)$ belongs to the cone of positive
semidefinite matrices of size $s(r)$, which is closed.
Next, and again by pointwise convergence, for every $s\in\N$,
\[\sum_{j=0}^s \sum_{i=1}^n \,(y^*)^{(i)}_{2j}/j{\rm !}\,=\,
\lim_{k\to\infty}\:\sum_{j=0}^s\sum_{i=1}^n \,(y^{(r_k)})^{(i)}_{2j}/j{\rm
!}\,\leq\,n\e ^{M^2},\]
and so, by the Monotone Convergence Theorem
\begin{equation}
\label{w11}
\sum_{j=0}^\infty\sum_{i=1}^n \,(y^*)^{(i)}_{2j}/j{\rm !}
\,=\,\lim_{s\to\infty}\sum_{j=0}^s\sum_{i=1}^n \,(y^*)^{(i)}_{2j}/j{\rm !}
\,\leq\,
n\e ^{M^2}.
\end{equation}

But (\ref{w11}) implies that $\y^*$ satisfies Carleman's condition (\ref{suff-1}) is
satisfied. Indeed, from (\ref{w11}), for all $i=1,\ldots,n$, we have
$(y^*)^{(i)}_{2k}<\rho k{\rm !}$ for all $k\in\N$, and so, as
$k{\rm!}\leq k^k=\sqrt{k}^{2k}$,
\[[(y^*)^{(i)}_{2k}]^{-1/2k}>(\rho)^{-1/2k}/\sqrt{k},\]
which in turn implies 
\[\sum_{k=0}^\infty [(y^*)^{(i)}_{2k}]^{-1/2k}>\sum_{k=0}^\infty \frac{\rho^{-1/2k}}{\sqrt{k}}\,=\,+\infty.\]
Hence, by Theorem \ref{suff}, $\y^*$ is a determinate moment sequence,
that is, there exists a unique measure $\mu ^*$ on $\R ^n$, such that
\[y^*_\alpha\,=\,\int_{\R^n}x^\alpha\,d\mu ^*,\qquad \alpha\in\N^n.\]
By (\ref{w11}),
\[\int \,\sum_{i=1}^n \e ^{x_i ^2}\,d\mu ^*\,=\,
\sum_{j=0}^\infty\sum_{i=1}^n \,(y^*)^{(i)}_{2j}/j{\rm !}\,\leq\,n\e ^{M^2},\]
which proves that $\mu ^*$ is admissible for $\mathcal{P}_M$.

But then, again by the pointwise convergence
(\ref{pointwise}) of $\y^{(r_k)}$ to $\y ^*$, we get
$L_{\y^{(r_k)}}(f)\to L_{\y^*}(f)=\int fd\mu ^*$ 
as $k\to\infty$, which, in view of
$L_{\y^{(r_k)}}(f)\leq \inf\mathcal{P}_M$ for all $k$, implies 
\[\int f\,d\mu ^*\,=\,L_{\y^*}(f)\,\leq \,\inf\mathcal{P}_M.\]
But this proves that $\mu ^*$ is an optimal solution of
$\mathcal{P}_M$ because $\mu ^*$ is admissible for $\mathcal{P}_M$ with
value $\int f d\mu ^*\leq \inf\mathcal{P}_M$. 
As the converging subsequence $\{r_k\}$
was arbitrary, it is true for every limit point.
Hence, we have proved (i) and (ii).

(iii) Let $\y$ be the moment sequence associated with the
probability measure $\mu$ on the ball
\[B_{M/2}\,=\,\{\y=\{y_\alpha\}\in l_\infty\,\vert\quad \Vert
\y\Vert_\infty\leq M/2\},\]
with uniform distribution. That is,
\[\mu(B)\,=\,M^{-n}\,\int_{B\cap B_{M/2}}\,dx,\qquad
B\in\mathcal{B},\]
where $\mathcal{B}$ is the sigma-algebra of Borel subsets of $\R^n$.

As $\mu$ has a continuous density $f_\mu>0$ 
on $B_{M/2}$, it follows easily that $M_r(\y)\succ0$
for all $r\geq r_f$. In addition, 
\[\sum_{k=0}^ r\sum_{i=1}^n y^{(i)}_{2k}/k{\rm !}\,<\,
\int \sum_{i=1}^n \e ^{x^2_i}\,d\mu\,<\,n\e
^{M^2},\]
so that $\y$ is a strictly admissible solution for $\Q_r$. Hence, the
SDP problem $\Q_r$ satisfies Slater's condition,
and so, there is no duality gap between
$\Q_r$ and $\Q^*_r$, and $\Q^*_r$ is solvable 
if $\inf\Q_r$ is finite; see e.g. Vandenberghe and Boyd \cite{boyd}.
Thus, $\Q^*_r$ is solvable because we proved that $\Q_r$ is solvable. 
In  other words, $\sup\Q^*_r=\max\Q^*_r=\min\Q_r$, the
desired result.
$\qed$

\subsection{Proof of Theorem \ref{uncons}}
\label{proofuncons}

It suffices to prove (i) and (ii) for the case $f ^*>0$.
Indeed, if $f^*=0$ take $\epsilon>0$ arbitrary, fixed.
Then $f+n\epsilon\geq f^*_\epsilon=f^*+n\epsilon>0$ and so, suppose that (\ref{uncons-1}) holds for 
$f+n\epsilon$ (for some $r_0,\lambda_0$). 
In particular, pick $\lambda\geq\lambda_0+\epsilon$,
so that
\[f+n\epsilon+(\lambda-\epsilon)
\sum_{k=0}^{r_{\lambda}}\sum_{j=1}^n
\frac{x_j^{2k}}{\k}\,=\,q_\lambda,\]
(with $q_\lambda$ s.o.s.), for $r_{\lambda}\geq r_0$. 
Equivalently,
\[f+\lambda\sum_{k=0}^{r_{\lambda}}\sum_{j=1}^n
\frac{x_j^{2k}}{\k}\,=\,q_\lambda+
\epsilon\sum_{k=1}^{r_{\lambda}}\sum_{j=1}^n
\frac{x_j^{2k}}{\k}\,=\,\hat{q}_\lambda,\]
where $\hat{q}_\lambda$ is a s.o.s. Hence (\ref{uncons-1}) also holds
for $f$ (with $\lambda_0+\epsilon$ in lieu of $\lambda_0$). 

Similarly, for (\ref{uncons-2}). As $f^*=0$, 
$f+n\epsilon>0$ and so, suppose that (\ref{uncons-2}) holds for 
$f+n\epsilon$. In particular,
\[f+n\epsilon +
\epsilon\sum_{k=0}^{r_{\epsilon}}\sum_{j=1}^n
\frac{x_j^{2k}}{\k}\,=\,q_\epsilon,\]
(with $q_\epsilon$ s.o.s.), for some $r_{\epsilon}$. 
Equivalently,
\[f+2\epsilon\sum_{k=0}^{r_{\epsilon}}\sum_{j=1}^n
\frac{x_j^{2k}}{\k}\,=\,q_\epsilon+
\epsilon\sum_{k=1}^{r_{\epsilon}}\sum_{j=1}^n
\frac{x_j^{2k}}{\k}\,=\,\hat{q}_\epsilon,\]
where $\hat{q}_\epsilon$ is a s.o.s. Hence (\ref{uncons-2}) also holds
for $f$. Therefore, we will assume that $f^*>0$.

(i) As $f^*>0$, let $M_0$ be such that $f^*>1/M_0$, and fix $M>M_0$.
Consider the SDP problem $\Q^*_r$ defined in (\ref{dualagain}),
associated with $M$. By Proposition \ref{prop1}, $f^*\leq\inf\mathcal{P}_{M}$.
By Theorem \ref{th1}, $\max\Q^*_r=\min\Q_r\uparrow 
\inf\mathcal{P}_{M}\geq f^*$. Therefore, 
there exists some $r_M\geq r_f$ such that
$\max\Q^*_{r_M}\geq f^*-1/M>0$. That is, if 
$(q_M,\lambda_M,\gamma_M)$ is an
optimal solution of $\Q^*_{r_M}$, then 
$\gamma_M-n\lambda_M\e ^{M^2}\geq f^*-1/M>0$.
In addition, 
\[f-\gamma_M= q_M-\lambda_M
\sum_{k=0}^{r_M}\sum_{j=1}^n \frac{x_j^{2k}}{\k},\]
that we rewrite
\begin{equation}
\label{n1}
f-(\gamma_M-n\lambda_M\e ^{M^2})= 
q_M+\lambda_M\left(n\e ^{M^2}-
\sum_{k=0}^{r_M}\sum_{j=1}^n \frac{x_j^{2k}}{\k}\right).
\end{equation}
Equivalently,
\[f+\lambda_M\sum_{k=0}^{r_M}\sum_{j=1}^n \frac{x_j^{2k}}{\k}
\,=\, q_M+n\lambda_M\e ^{M^2} +(\gamma_M-n\lambda_M\e^{ M^2}).\]
Define $\hat{q}_M$ to be the s.o.s. polynomial
\[\hat{q}_M\,:=\,q_M+n\lambda_M\e ^{M^2}
+(\gamma_M-n\lambda_M\e^{ M^2}),\]
so that we obtain

\begin{equation}
\label{n2}
f+\lambda_M\sum_{k=0}^{r_M}\sum_{j=1}^n
\frac{x_j^{2k}}{\k}\,=\,\hat{q}_M,
\end{equation}
the desired result.

If we now take $r>r_M$ and $\lambda\geq \lambda_M$
we also have
\begin{eqnarray}
\nonumber
f+\lambda\sum_{k=0}^r\sum_{j=1}^n
\frac{x_j^{2k}}{\k}&=&
f+\lambda_M
\sum_{k=0}^{r_M}\sum_{j=1}^n\frac{x_j^{2k}}{\k}\\
\nonumber
&+&\lambda_M \sum_{k=r_M+1}^{r}\sum_{j=1}^n\frac{x_j^{2k}}{\k}+
(\lambda-\lambda_M) \sum_{k=0}^r\sum_{j=1}^n\frac{x_j^{2k}}{\k}\\
\nonumber
&=&\hat{q}_M+\lambda_M \sum_{k=r_M+1}^{r}\sum_{j=1}^n\frac{x_j^{2k}}{\k}
+(\lambda-\lambda_M) \sum_{k=0}^r\sum_{j=1}^n\frac{x_j^{2k}}{\k}\\
\nonumber
&=&\hat{\hat{q}}_M,
\end{eqnarray}
that is,
\begin{equation}
\label{n3}
f+\lambda\sum_{k=0}^r\sum_{j=1}^n
\frac{x_j^{2k}}{\k}\,=\,\hat{\hat{q}}_M,
\end{equation}
where $\hat{\hat{q}}_M$ is a s.o.s. polynomial, the desired result. $\qed$

(ii) Let $M$ be as in (i) above.
%Let $z\in B_{1}$, fixed, arbitrary and 
Evaluating (\ref{n1}) at $x=0$, and writing $f(0)=f(0)-f^*+f^*$, yields
\[f(0)-f^*+f^*-(\gamma_M-n\lambda_M\e ^{M^2})
\,=\,q_M(0)+
n\lambda_M\,(\e ^{M^2}-1),\]
%\sum_{k=0}^{r_M}\sum_{i=1}^n \frac{z_i^{2k}}{\k}\right),\]
and as $1/M>f^*-(\gamma_M-n\lambda_M\e ^{M^2})$,
\[\lambda_M\,\leq\,\frac{1/M+ f(0)-f^*}{n\,(\e ^{M^2}-1)}.\]
%\sum_{k=0}^{r_M}\sum_{i=1}^n\frac{z_i^{2k}}{\k}\right).\]
%Recall that $z\in B_{1}$. 
%Hence, with $0<\rho_M:=n(\e ^{M^2}-\e )\leq n\e ^{M^2}-
%\sum_{k=0}^{r_M}\sum_{i=1}^n \frac{z_i^{2k}}{\k},\]
%one obtains $0\leq \lambda_M\,\leq\,(1/M+
%f(z)-f^*)/\rho_M$. 
Now, letting $M\to\infty$, yields
$\lambda_M\to 0$.

Now, let $\epsilon >0$ be fixed, arbitrary. There is some
$M>M_0$ such that $\lambda_M\leq\epsilon$ in
(\ref{n2}). Therefore, (\ref{uncons-2}) is just
(\ref{n3}) with $\lambda:=\epsilon>\lambda_M$ and
$r=r_\epsilon\geq r_M$.
Finally, from this, we immediately have
\[\Vert f-f_\epsilon\Vert_1\,\leq\,\epsilon
\sum_{i=1}^n\sum_{k=0}^\infty \frac{1}{\k}\,=\,\epsilon \,n \e \to
0,\quad \mbox{as }\epsilon\downarrow 0.\]
$\qed$

\section{Appendix}
\label{appendix}
In this section we derive two auxiliary results that are helpful in
the proofs of Theorem \ref{th1} and Theorem \ref{uncons} in \S\ref{proofs}.

\begin{lemma}
\label{lemma2}
Let $n=2$ and let $\y$ be a sequence indexed in the basis 
(\ref{aa4}), and such that $M_r(\y)\succeq 0$. Then all the diagonal entries of $M_r(\y)$
are bounded by $\tau_r:=\max_{k=1,\ldots,r}\max\,[y_{2k,0},y_{0,2k}]$.
\end{lemma}
\begin{proof}
It suffices to prove that all the entries $y_{2\alpha,2\beta}$ with
$\vert\alpha+\beta\vert=r$ are bounded by
$s_r:=\max\,[y_{2r,0},y_{0,2r}]$, and repeat the argument for 
entries $y_{2\alpha,2\beta}$ with $\vert\alpha+\beta\vert=r-1,r-2,$
etc ... Then, take $\tau_r:=\max_{k=1,\ldots,r}s_k$. 
So, consider the odd case $r=2p+1$, and the even case $r=2p$.

- The odd case $r=2p+1$. 
Let $\Gamma:=\{(2\alpha,2\beta)\,\vert\quad\alpha+\beta
=r,\:\alpha\neq 0\}$, and notice that 
\[\Gamma\,=\,\{(2r-2k,2k)\:\vert \quad k=1,\ldots,r-1\}\,=\,\Gamma_1\cup\Gamma_2\]
with
\[\Gamma_1\,:=\,\{(r,0)+(r-2k,2k),\vert \quad k=1,\ldots,p\},\]
and
\[\Gamma_2:=\{(0,r)+(2j,r-2j),\vert\quad j=1,\ldots,p\}.\]
Therefore, consider the two rows (and columns) corresponding to
the indices $(r,0)$ and $(r-2k,2k)$, or $(0,r)$ and $(2j,r-2j)$.
In view of $M_r(\y)\succeq 0$, one has
\begin{equation}
\label{enfin}
\left\{\begin{array}{lcl}
y_{2r,0}\times y_{2r-4k,4k}&\geq &(y_{2r-2k,2k})^2,\qquad
k=1,\ldots,p,\\
y_{0,2r}\times y_{4j,2r-4j}&\geq &(y_{2j,2r-2j})^2,\qquad
j=1,\ldots,p.
\end{array}\right.
\end{equation}
Thus, let $s:=\max\,\{y_{2\alpha,2\beta}\,\vert \quad
\alpha+\beta=r,\alpha\neq 0\}$, so that either
$s=y_{2r-2k^*,2k^*}$ for some $1\leq k^*\leq p$, or
$s=y_{2j^*,2r-2j^*}$ for some $1\leq j^*\leq p$.
But then, in view of (\ref{enfin}), and with
$s_r:=\max\,[y_{2r,0},y_{0,2r}]$,
\[s_r \times s\geq  y_{2r,0}\times y_{2r-4k^*,4k^*}\geq
(y_{2r-2k^*,2k^*})^2=s ^2,\]
or,
\[s_r \times s\geq  y_{0,2r}\times y_{4j^*,2r-4j^*}\geq
(y_{2j^*,2r-2j^*})^2=s ^2,\]
so that $s\leq s_r$, the desired result.

- The even case $r=2p$.
Again, the set $\Gamma:=\{(2\alpha,2\beta)\,\vert\quad\alpha+\beta
=r,\:\alpha\neq 0\}$ can be written $\Gamma=\Gamma_1\cup\Gamma_2$, with
\[\Gamma_1\,:=\,\{(r,0)+(r-2k,2k),\vert \quad k=1,\ldots,p\},\]
and
\[\Gamma_2:=\{(0,r)+(2j,r-2j),\vert\quad j=1,\ldots,p\}.\]
The only difference with the odd case is that
$\Gamma_1\cap\Gamma_2=(p,p)\neq\emptyset$. But
the rest of the proof is the same as in the odd case.
\end{proof}

\begin{lemma}
\label{lemma1}
Let $r\in\N$ be fixed, and let $\y$ be a sequence such that
the associated moment matrix $M_r(\y)$ is positive semidefinite, i.e.,
$M_r(\y)\succeq0$. Assume that there is some $\tau_r\in\R$ such that
the diagonal elements $\{y ^{(i)}_{2k}\}$ satisfy
$y ^{(i)}_{2k}\leq\tau_r$, for all $k=1,\ldots,r$, and all $i=1,\ldots,n$.

Then, the diagonal elements of $M_r(\y)$ are all bounded by $\tau_r$
(i.e., $y_{2\alpha}\leq \tau_r$ for all $\alpha\in\N^n$, with $\vert\alpha\vert\leq r$).
\end{lemma}

\begin{proof}
The proof is by induction on the the number $n$ of variables.
By our assumption it is true for $n=1$, and by Lemma \ref{lemma2}, it is
true for $n=2$. Thus, suppose it is true for $k=1,2,\ldots,n-1$ variables and
consider the case of $n$ variables (with $n>3$).

By our induction hypothesis, it is true for all elements $y_{2\alpha}$ where at least
one index, say $\alpha_i$, is zero ($\alpha_i=0$).
Indeed, the submatrix $A ^{(i)}_r(\y)$ of $A_r(\y)$, obtained 
from $A_r(\y)$ by deleting all rows and columns corresponding to indices
$\alpha\in\N^n$ in the basis (\ref{aa4}),
with $\alpha_i>0$, is a moment matrix of order $r$, with
$n-1$ variables $x_1,\ldots,x_{i-1},x_{i+1},\ldots,x_n$.  Hence, by
a permutation of rows and columns, we can write
\[A_r(\y)\,=\,\left[\begin{array}{ccc}
A ^{(i)}_r(\y)&\vert& B\\
-&\vert &-\\
B' &\vert & C\end{array}\right],\]
for some appropriate matrices $B$ and $C$.
In particular, all elements $y_{2\alpha}$ with $\alpha_i=0$,
are diagonal elements of $A ^{(i)}_r(\y)$.
In addition, its diagonal elements $y^{(j)}_{2k}$, $j\neq i$, 
are all bounded by $\tau_r$. And of course, $A ^{(i)}_r(\y)\succeq 0$.
Therefore, by our induction hypothesis, all its diagonal elements 
are bounded by $\tau_r$. As $i$ was arbitrary, we conclude that all
elements $y_{2\alpha}$ with at least one index being zero, are all
bounded by $\tau_r$.

We next prove it is true for an
arbitrary element $y_{2\alpha}$ with $\vert\alpha\vert\leq r$ and
$\alpha>0$, i.e., $\alpha_j\geq 1$ for all $j=1,\ldots,n$.
With no loss of generality, we assume that $\alpha_1\leq\alpha_2\leq
\ldots\leq\alpha_n$.

Consider the two elements $y_{2\alpha_1,0,\beta}$ and
$y_{0,2\alpha_2,\gamma}$, with
$\beta,\gamma\in\N^{n-2}$ such that:
\[\vert\beta\vert\,=\,\vert\alpha\vert-2\alpha_1;\quad
\vert\gamma\vert\,\,=\,\vert\alpha\vert-2\alpha_2,\]
and
\[(2\alpha_1,0,\beta)+(0,2\alpha_2,\gamma)\,=\,(2\alpha_1,2\alpha_2,\beta+\gamma)\,=\,
2\alpha.\]
So, for instance, take $\beta=(\beta_3,\beta_4,\ldots,\beta_n)$, 
$\gamma=(\gamma_3,\gamma_4,\ldots,\gamma_n)$, 
defined by
\[\beta\,:=\,(\alpha_3+\alpha_2-\alpha_1,\alpha_4,\ldots,\alpha_n),\quad
\gamma\,:=\,(\alpha_3+\alpha_1-\alpha_2,\alpha_4,\ldots,\alpha_n).\]

By construction, we have 
$4\alpha_1+2\vert\beta\vert= 4\alpha_2+2\vert\gamma\vert=2\vert\alpha\vert\leq 2r$,
so that both
$y_{4\alpha_1,0,2\beta}$ and 
$y_{0,4\alpha_2,2\gamma}$ are diagonal elements of $M_r(\y)$ with at
least one entry equal to $0$. Hence, by the induction hypothesis,
\[y_{4\alpha_1,0,2\beta}\,\leq\,\tau_r,\quad
y_{0,4\alpha_2,2\gamma}\,\leq\,\tau_r.\]
Next, consider the two rows and columns indexed by 
$(2\alpha_1,0,\beta)$ and $(0,2\alpha_2,\gamma)$.
The constraint $M_r(\y)\succeq 0$ clearly implies
\[\tau_r^2\,\geq\,y_{4\alpha_1,0,2\beta}\times y_{0,4\alpha_2,2\gamma}\,\geq\,
(y_{2\alpha_1,2\alpha_2,\beta+\gamma})^2\,=\,y_{2\alpha}^2.\]
Hence, $y_{2\alpha}\leq \tau_r$, the desired result.
\end{proof}

\bibliographystyle{amsplain}

\end{document}